\documentclass[
preprintnumbers,amsmath,amssymb]{revtex4}

\usepackage{graphicx}
\usepackage{dcolumn}
\usepackage{bm}

\usepackage{graphicx}
\usepackage{epsfig}
 \usepackage{color}
\begin{document}

\title{\Large Integrable Abel equations and Vein's Abel equation}

\author{Stefan C. Mancas}
\email{stefan.mancas@erau.edu}
\affiliation{%
Fakult\"{a}t f\"{u}r Informatik und Mathematik\\
Hochschule M\"{u}nchen - Munich University of Applied Sciences, Germany
}%

\author{Haret C. Rosu}
\email{hcr@ipicyt.edu.mx}
\affiliation{IPICyT, Instituto Potosino de Investigacion Cientifica y Tecnologica,\\
Camino a la presa San Jos\'e 2055, Col. Lomas 4a Secci\'on, 78216 San Luis Potos\'{\i}, S.L.P., Mexico}

\date{online at Math. Meth. Appl. Sci. -- 7/28/2015}

\begin{center}
Math. Meth. Appl. Sci. 39 (2016) 1376-1387
\end{center}


\begin{abstract}
\noindent We first reformulate and expand with several novel findings some of the basic results in the integrability of Abel equations.
Next, these results are applied to Vein's Abel equation whose solutions are expressed in terms of the third order hyperbolic functions and a phase space analysis of the corresponding nonlinear oscillator is also provided.\\

\end{abstract}

\maketitle

\section{Introduction}\label{secI}

The first order Abel nonlinear equation of the second kind has the form
  \begin{equation}\label{d-3}
(w+s)\frac{dw}{dx}+p+q_1w+q_2w^2+rw^3=0~,
 \end{equation}
where $p$, $q_1$, $q_1$, $r$, and $s$ are all some functions of $x$. The case $r=0$ has been introduced more than two hundred years ago by Abel \cite{a1}.
Many solvable equations of this type are collected in \cite{polzai} and other ones can be found in more recent works \cite{ct1,ct2,mak1,mak2,pz,sh,mr}. The transformation $1/y =w+s$ converts this equation in the form
  \begin{equation}\label{d-4}
 \frac{dy}{dx}=r+(q_2-3rs)y+(q_1-\frac {ds}{dx}-2q_2s+3rs^2)y^2+(p-q_1s+q_2s^2-rs^3)y^3~,
  \end{equation}
which is Abel's equation of the first kind. We see that Abel's original equation of the second kind is actually a homogeneous case of that of the first kind.

\medskip

One can also write \eqref{d-4} in the canonical form
\begin{equation}\label{n2p}
\frac{dy}{dx}=f_0(x)+f_1(x)y+f_2(x)y^2+f_3(x)y^3=F(x,y)~,
\end{equation}
where only the case $s=0$ leads to the simple identifications: $f_0=r$, $f_1=q_2$, $f_2=q_1$, and $f_3=p$.

The integrability features of Abel's equation in the canonical form are extremely important because of its connection with nonlinear second order differential equations that phenomenologically describe a wide class of nonlinear oscillators which in this way can be treated analytically.
In particular, Vein \cite{ve} introduced an Abel equation with $f_0=0$ and rational forms of $f_1$, $f_2$, and $f_3$ whose solutions are expressed in terms of third order hyperbolic functions. However, to the best of our knowledge, the paper of Vein went unnoticed for many years and only recently Yamaleev \cite{yama1,yama2} provided a generalization in the framework of third order multicomplex algebra. Needless to say, the corresponding Vein\rq{}s nonlinear oscillator has not been studied in the literature. This was the main motivation for writing this paper, which is organized as follows. In section \ref{sec1}, we present a lemma that encodes the connection of Abel equation with nonlinear oscillator equations and some related simple results. In section \ref{sec2}, we reformulate several integrability results for Abel's equation with the purpose to apply them to Vein's Abel equation, which is the subject of section \ref{sec3}. The dynamical systems analysis of Vein's oscillator is developed in section \ref{sec4}, and we end up with the conclusions.

\section{Connection with second order nonlinear ODEs}\label{sec1}

The known importance of Abel's equation in its canonical form \eqref{n2p} stems from the fact that its integrability leads to closed form solutions to a general nonlinear ODE of the form
\begin{equation}\label{n1p}
x_{\zeta \zeta}+f_2(x) x_{\zeta}  +f_3(x)+f_1(x) x^2_{\zeta}+f_0(x) x^3_{\zeta}=0~,
\end{equation}
where the variables $x(\zeta)$ and $y(x(\zeta))$ are some parametric solutions that depend on a
generalized coordinate $\zeta$.
This can be expressed by the following lemma \cite{mr}

\begin{quote}
{\bf Lemma 1}. Solutions to a general second-order ODE of type (\ref{n1p}) may be obtained via the solutions to Abel's equation (\ref{n2p}) and vice versa using the following relationship
%
\begin{equation}\label{n10}
\frac{dx}{d\zeta}=v(x(\zeta))~.
\end{equation}

\end{quote}

{\bf Proof}. To show the equivalence, one just needs the chain rule
\begin{equation}\label{eq1}
\frac{d^2x}{d\zeta^2}=\frac{d v}{d x}\frac{dx}{d\zeta}=v \frac{d v}{d x}
\end{equation}
which turns (\ref{n1p}) into the Abel equation of the second kind in canonical form
\begin{equation}\label{eq2}
v \frac{d v}{dx}+f_3(x)+f_2(x)v+f_1(x)v^2+f_0(x)v^3=0.
\end{equation}
Via the inverse transformation
\begin{equation}\label{n10bis}
v(x(\zeta))=\frac{1}{y(x(\zeta))}~.
\end{equation}
of the dependent variable,
equation (\ref{eq2}) becomes (\ref{n2p}).
Moreover, the linear term in (\ref{n2p}) can always be eliminated via the transformation
\begin{equation}\label{tr1}
y(x)=z(x)e^{\int f_1(x) dx}
\end{equation}
which gives
\begin{equation}\label{n2pp}
\frac{dz}{dx}=h_0(x)+h_2(x)z^2+h_3(x)z^3=H(x,z)~,
\end{equation}
where
\begin{equation}\label{pee}
\begin{aligned}
	h_0(x)&=f_0(x)e^{-\int f_1(x)dx}\\
	 h_2(x) &=f_2(x)e^{\int f_1(x)dx}\\
	  h_3(x) &=f_3(x)e^{2\int f_1(x)dx}~.\\
	 \end{aligned}
	 \end{equation}

The case where $f_0=0 \rightarrow r=0$ is seen from (\ref{d-4}) to be the one actually considered by Abel, for in this case the reduced Abel equation
  \begin{equation}\label{d}
(w+s)\frac{dw}{dx}+p+q_1w+q_2w^2=0~,
 \end{equation}
can be always put in the form
 \begin{equation}\label{n2ppp}
\frac{dz}{dx}=h_2(x)z^2+h_3(x)z^3~
\end{equation}
where $h_2(x)$ and $h_3(x)$ are reduced functions and are given by the expressions
\begin{equation}\label{pee2}
\begin{aligned}
	 h_2(x) &=(q_1-\frac{ds}{dx}-2q_2s)e^{\int q_2(x)dx}\\
	  h_3(x) &=(p-q_1s+q_2s^2)e^{2\int q_2(x)dx}~.\\
	 \end{aligned}
	 \end{equation}
Using the lemma, the reduced Abel equation corresponds to a linear ODE without higher-order dissipative terms
\begin{equation}\label{m1p}
x_{\zeta \zeta}+h_2(x)x_{\zeta}  +h_3(x)=0~.
\end{equation}

    On the other hand, in the case of $f_3=0$, the Abel equation of the first-kind (\ref{n2p}) becomes a Riccati equation, while the second-kind homogeneous Abel equation \eqref{eq2} is reduced to the Riccati equation
\begin{equation}\label{eqabric}
\frac{d v}{dx}=-f_2(x)-f_1(x)v-f_0(x)v^2~,
\end{equation}
which is equivalent to
\begin{equation}\label{n2r}
\frac{dy}{dx}=f_0(x)+f_1(x)y+f_2(x)y^2~.
\end{equation}

This corresponds to the simple fact that an inverse power of a Riccati solution also satisfies a Riccati equation with redistributed coefficients. Now, we eliminate the linear part to obtain the reduced Riccati equation
 \begin{equation}\label{n2pppp}
\frac{dz}{dx}=h_0(x)+h_2(x)z^2~
\end{equation}
which corresponds to a nonlinear ODE with higher-order dissipative terms
\begin{equation}\label{p1p}
x_{\zeta \zeta}+h_2(x)x_{\zeta}+h_0(x) x_{\zeta}^3=0~,
\end{equation}
where now the coefficients are
\begin{equation}\label{gee}
\begin{aligned}
	 h_0(x) &=re^{-\int (q_2- 3 r s )dx}\\
	  h_2(x) &=(q_1-s\rq{}-2 q_2s+3rs^2)e^{\int (q_2- 3 r s )dx}~.\\
	 \end{aligned}
	 \end{equation}

\section{Some Abel integrability cases}\label{sec2}

\subsection{Abel's equation of constant coefficients}

Denote the coefficients $f_i(x)=A_i$ of (\ref{n2p}), where $A_i \in \mathbb{R}$ are constants and $A_3 \ne 0$ so that $F(x,y)=F(y)$. It is obvious that the roots of the equation $F(y)=0$ are themselves solutions of (\ref{n2p}). More generally, the general solution of (\ref{n2p}) is obtained via factorization of the denominator in the right-hand side of 
\begin{equation}\label{d-5a}
\frac{dy}{y^3+\frac{A_2}{A_3}y^2+\frac{A_1}{A_3} y+\frac{A_0}{A_3}}=A_3dx~,
\end{equation}
which leads to the following cases:
\begin{itemize}
\item [(i)]  If $y_1 \ne y_2 \ne y_3 \in \mathbb{R} $ , 
so integration of (\ref{d-5a}) has the form

\begin{equation}\label{d-5b}
\int{\frac{dy}{(y-y_1)(y-y_2)(y-y_3)}}= (y-y_1)^{y_2-y_3}(y-y_2)^{y_3-y_1}(y-y_3)^{y_1-y_2}=c_1e^{A_3x}~,
\end{equation}

\item[(ii)] If $y_1 \ne y_2= y_3 \in \mathbb{R}$, 
then the integration of (\ref{d-5a}) is of the form

\begin{equation}\label{d-5c}
\int{\frac{dy}{(y-y_1)(y-y_2)^2}}= \frac{1}{(y_1-y_2)(y-y_2)}+\frac{1}{(y_1-y_2)^2}\ln{|\frac{y-y_1}{y-y_2}|}=A_3x +c_2~,
\end{equation}

\item [(iii)]  If $y_1 =y_2= y_3 \in \mathbb{R}$ , 
so (\ref{d-5a}) becomes

\begin{equation}\label{d-5d}
\int{\frac{dy}{(y-y_1)^3}}= \frac{1}{-2(y-y_1)^2}=A_3x +c_3~,
\end{equation}

\item [(iv)] If $y_1= \overline{{y_2}}=\alpha+ i \beta \in \mathbb{C}$ and $y_3 \in \mathbb{R}$, 
the integration of (\ref{d-5a}) leads to

\begin{equation}\label{d-5e}
\int{\frac{dy}{\Big((y-\alpha)^2+\beta^2\Big)(y-y_3)}}= \ln\Big|\sqrt{(y-\alpha)^2+\beta^2}\Big|+\frac{\alpha-y_1}{\beta}\arctan\frac{y-\alpha}{\beta}=A_3x +c_4~.
\end{equation}

\end{itemize}

 \subsection{Integrability based on the normal form of Abel's equation}
If the following transformations as given in Kamke's book \cite{kam}
\begin{equation}\label{tr}
\begin{aligned}
     y(x)&=\omega(x)\eta(\xi(x))-\frac{f_2(x)}{3 f_3(x)}~, \qquad \omega(x)=e^{\int{\left(f_1-\frac{f_2^2}{3 f_3}\right)dx}}~,\\
	 \xi(x) &= \int {f_3\,\omega^2\,dx}\\
\end{aligned}
\end{equation}
are applied to equation \eqref{n2p}, then one obtains Abel\rq{}s equation in normal form
 \begin{equation}\label{nf}
\frac{d \eta}{d \xi}=\eta^3+I(x)~,
\end{equation}
where the invariant $I(x)$ is given by
\begin{equation}\label{in}
I(x)=\frac {f_0+\frac 1 3\frac {d}{dx}\left(\frac{f_2}{f_3}\right)-\frac{f_1f_2}{3 f_3}+\frac{2 f_2^3}{27 f_3^2}}{f_3 \omega^3}~.
\end{equation}
Thus, we conclude that if $I(x) \equiv {\rm constant}$, then \eqref{nf} is integrable since it is separable.
If one chooses relations between the functions $f_i$ such that the invariant is null and letting $f_2(x)/f_3(x)=n(x)$, then
\begin{equation}
\frac{dn}{dx}+\frac 2 9 f_2 n^2-f_1n+3 f_0=0
\end{equation}
is a Riccati equation, which is always integrable because it is obtained from Abel\rq{}s equation in normal form, which is integrable.
Thus, the solution to \eqref{nf} is
\begin{equation}
\eta=\frac{1}{\sqrt{c- 2 \xi}}~,
\end{equation}
and explicitly, the solution to \eqref{n2p} with null invariant $I(x)$ is
\begin{equation}
y=\frac{e^{  \int{\left (f_1-\frac{f_2^2}{3 f_3}\right)dx}}}{\sqrt{c- 2 \displaystyle\int  f_3 \left(e^{\int{(f_1-\frac{f_2^2}{3 f_3})dx}}\right)^2}}-\frac {f_2}{3 f_3}~.
\end{equation}

\subsection{Integrability of Abel equation with non-constant coefficients and $f_0=0$}
In this subsection, we will consider the original Abel equation (\ref{d-4}) of the first-kind
\begin{equation}\label{Abel}
\frac{dy}{dx}=f_1(x)y+f_2(x)y^2+f_3(x)y^3~,
\end{equation}
which corresponds to (\ref{n1p}) without the cubic nonlinearity.

Let us use the transformation
\begin{equation}
y(x)=\frac{e^{\int f_1(x)dx}}{\nu(x)}~,
\end{equation}
which is \eqref{tr1} with $z=\frac{1}{\nu}$
 that allows us to put \eqref{Abel} into a differential form
\begin{equation}\label{d-8}
\nu d \nu+(P+Q\nu)dx=0~,
\end{equation}
where
\begin{equation}\label{pe}
\begin{aligned}
	P(x)&=f_3(x)e^{2\int f_1(x)dx}\\
	 Q(x) &=f_2(x)e^{\int f_1(x)dx}~.
	 \end{aligned}
	 \end{equation}
Thus, the original equation \eqref{d-3} is considerably reduced only by having $f_0=0$.

\medskip

For equation \eqref{d-8}, if there exists an integrating factor \cite{davis}
\begin{equation}\label{d-9}
\mu(x,\nu)=e^{-k(\nu+\int Q(x)dx)}~,
\end{equation}
where $k$ is a constant, then
\begin{equation}\label{d-9a}
(P+Q\nu)\mu dx+\nu \mu d\nu=0~,
\end{equation}
is integrable, provided that
\begin{equation}\label{d-9b}
Q \mu +(P+Q\nu)\frac{\partial \mu}{\partial \nu}=  \nu \frac{\partial \mu}{\partial x}~,
\end{equation}
which leads to $Q=kP$. Hence, the necessary condition for the integration of (\ref{Abel}) is
 \begin{equation}\label{con}
 kf_3(x)e^{\int f_1 (x)dx}=f_2(x) \rightarrow \frac{d}{dx}\big(\frac{f_2}{f_3}\big) =f_1\frac{f_2}{f_3} \rightarrow \frac{dn}{dx}=f_1 n
 \end{equation}

Using $\mu(x,\nu)=e^{-k \nu}e^{-k^2\int {P(x) dx}}$, we find that there exists the potential $\Psi(x,\nu)=c_{\Psi}$, which satisfies
\begin{equation}\label{d-9c}
\begin{aligned}
	\frac{\partial \Psi}{\partial x} &= (1+k \nu)P \mu~, \\
	 \frac{\partial \Psi}{\partial \nu} &= \nu \mu.\\
\end{aligned}
\end{equation}
Therefore, the solution to (\ref{d-8}) is
\begin{equation}\label{sol2}
\Psi(x,\nu)=(1+k \nu)e^{-k(\nu+k\int Pdx)}=c_{\Psi},
\end{equation}
or in terms of the integrating factor
\begin{equation}\label{sol2a}
\Psi(x,\nu)=(1+k \nu)\mu(x,\nu)=c_{\Psi}.
\end{equation}
By substituting the condition \eqref{con}, the invariant \eqref{in} has the form
\begin{equation}
I(x)=\frac{2k^3}{27}e^{\int \frac{f_2^2}{f_3}dx}~.
\end{equation}

\subsection{Canonical form of Abel's equation and the integrating factor}
According to the book of Kamke \cite{kam}, for equations of the type \eqref{Abel} for which there is no condition (\ref{con}), one should change the variables according to
\begin{equation}\label{lem}
\begin{aligned}
	y(x)&= \tilde{\omega}(x)\tilde{\eta}(\zeta(x))~,\qquad \tilde{\omega}= e^{\int f_1(x) dx}~,\\
	 \zeta(x) &= \int {f_2(x)\tilde{\omega}(x)dx}~,\\
\end{aligned}
\end{equation}
which lead to the canonical form
\begin{equation}\label{norm1}
\frac{d \tilde{\eta}}{d \zeta}=\tilde{\eta}^2+g(\zeta)\tilde{\eta}^3,
\end{equation}
where
\begin{equation}
g(\zeta(x))=\frac{f_3(x)}{f_2(x)}e^{\int f_1(x) dx}\
\end{equation}
is the Appell invariant. Then, the integrating factor $\mu$ can be used to formulate the following interesting result:\\

{\bf Lemma 2}. Any Abel equation in the canonical form \eqref{norm1} is integrable as long as the invariant is constant, with $g(\zeta)= \frac 1k$, and solution given by \eqref{sol2}.

{\bf Proof}.  In \eqref{d-8}, let $\nu(x)=\frac{1}{\tilde{\eta}(\zeta(x))}$, where $\frac{d\zeta(x)}{dx}=Q(x)$, which leads to

\begin{equation}\label{norm2}
\frac{d \tilde{\eta}}{d \zeta}=\tilde{\eta}^2+\frac{P(x)}{Q(x)}\tilde{\eta}^3=\tilde{\eta}^2+ \frac{f_3(x)e^{\int f_1(x) dx}}{f_2(x)}\tilde{\eta}^3=\tilde{\eta}^2+\frac 1k\tilde{\eta}^3,
\end{equation}
that is separable, with solution
\begin{equation}\label{sol3}
\frac 1 k \ln\left |\frac{1}{\tilde{\eta}}+ \frac 1 k \right |=\zeta+c+\frac{1}{\tilde{\eta}}~,
\end{equation}
which after simplification leads back to \eqref{sol2}.

\section{Vein's Abel equation}\label{sec3}

The following Abel equation:
\begin{equation}\label{Abel2}
\frac{dy}{dx}=\frac{-2 b}{bx+a^2}y+\frac{3(ax+b^2)}{bx+a^2}y^2+\frac{x^3-3abx-a^3-b^3}{bx+a^2}y^3~,
\end{equation}
where $a, b \in \mathbb{R}$, 
is known from a paper of Vein to be integrable \cite{ve}.
We first notice that if we define the following determinants,
$$
D_3=\left| \begin{array}{rrr}
x & -a & -b\\-b & x & -a\\-a & -b & x
\end{array}\right|~, \qquad D_2=\left| \begin{array}{rr}
x & -b \\-b & -a\end{array}\right|~, \qquad D_1=|-b|~, \qquad D=\left| \begin{array}{rr}
-a & -b \\x & -a\end{array}\right|~,
$$
then they allow us to write the coefficients of Vein's Equation \eqref{Abel2} as $f_1=2D_1/D$, $f_2(x)=-3D_2/D$, and $f_3(x)=D_3/D$. Thus, \eqref{Abel2} can be written in the compact form
\begin{equation}\label{vdet}
\frac{dy}{dx}=\frac{2D_1}{D}y-\frac{3D_2}{D}y^2+\frac{D_3}{D}y^3~.
\end{equation}
Next, using the transformation given by \eqref{tr} with the coefficients from \eqref{Abel2}, we have
\begin{equation}\label{ome}
\omega=\frac{a^2+bx}{x^3-3 a b x-a^3-b^3}
\end{equation}
and
\begin{equation}\label{xii}
\frac{d \xi}{dx}=f_3(x)\omega^2=\omega.
\end{equation}
Using these expressions,  the invariant $I(x)$ becomes
\begin{equation}\label{in2}
I(x)=\frac {\frac {d}{dx}\frac{f_2}{3f_3}-\frac{f_1f_2}{3 f_3}+\frac{2 f_2^3}{27 f_3^2}}{f_3 \omega^3}=\frac{-\omega^2}{\frac {1}{\omega} \omega^3}\equiv -1~.
\end{equation}

Therefore, the normal  form of Abel's Equation \eqref{norm1} is
\begin{equation}\label{norm4}
\frac{d \eta}{d \xi}=\eta^3-1.
\end{equation}
Using \eqref{tr} and integrating \eqref{norm4}, we obtain
\begin{equation}\label{sol5}
\xi(\eta)-c=\frac 1 6 \ln \frac{(\eta-1)^2}{1+\eta+\eta^2}-\frac {\sqrt 3}{3}\arctan \frac {\sqrt 3}{3}(1+2 \eta)~.
\end{equation}
But because $\xi(x)=\int \omega  d x$, we find
\begin{equation}\label{xi}
\xi(x)-c=\frac 1 6  \ln {\frac{\big(x-(a+b)\big)^2}{x^2+(a+b)x+a^2-ab +b^2}}-\frac {\sqrt 3}{3}\arctan \frac {\sqrt 3}{3}{\frac{2x+a+b}{ a-b}}~.
\end{equation}
From the last two equations, one might infer that
\begin{equation}\label{rd}
\eta(x)=\frac{x+b}{a-b}
\end{equation}
 and consequently using \eqref{tr} would get
\begin{equation}\label{sd}
y(x)=\frac{b}{\big(x-(a+b)\big)(a-b)}~.
\end{equation}
 Unfortunately, because for \eqref{rd}, the argument of $\arctan$ is satisfied, but the argument of $\ln$ is not satisfied, then we conclude that \eqref{sd} is not the solution of \eqref{Abel2}.

\medskip

\begin{figure} [htb!]
 \centering
{\includegraphics{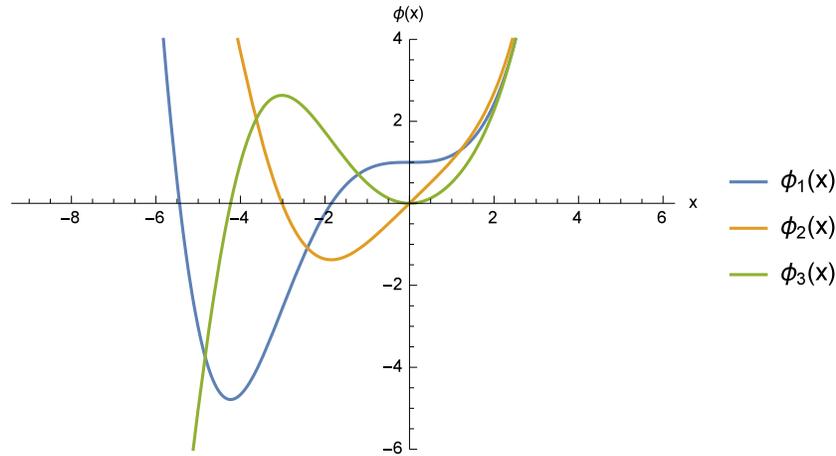}}
 \caption{The third-order hyperbolic functions $\phi_i$, $i=1,2,3$ of Equation~\eqref{triad}.}
 \label{f1}
\end{figure}

To see how the implicit solutions \eqref{sol5}, \eqref{rd} can be untangled and written explicitly, we will proceed as Vein did. Firstly, we  show  that there are actually three solutions of the type \eqref{sd}, which are the solutions of (\ref{Abel2}) and can be put in the form \cite{ve}:
\begin{equation}\label{ve-15}
\frac{b}{y_1}=(ax+b^2)-(bx+a^2)\bigg[\frac{t_3}{t_1}\bigg]~,
\qquad \frac{b}{y_2}=(ax+b^2)-(bx+a^2)\bigg[\frac{t_1}{t_2}\bigg]~,
\qquad \frac{b}{y_3}=(ax+b^2)-(bx+a^2)\bigg[\frac{t_2}{t_3}\bigg]~,
\end{equation}
with the $t$-functions defined cyclically as
$$
t_1=c\phi_1(s)+\phi_2(s)~,\qquad  t_2=c\phi_2(s)+\phi_3(s)~, \qquad t_3=c\phi_3(s)+\phi_1(s)~,
$$
where $c$ is an arbitrary constant as it does not appear in the differential equation. The $\phi$ functions are the
following triad:
\begin{align}\label{triad}
\phi_1(x)=&\frac{1}{3}\big[e^x+2e^{-x/2}\cos\frac{x\sqrt{3}}{2}\big]~, \nonumber\\
\phi_2(x)=&\frac{1}{3}\big[e^x-2e^{-x/2}\cos \big(\frac{x\sqrt{3}}{2}+\frac \pi 3\big)\big]~,\\
\phi_3(x)=&\frac{1}{3}\big[e^x-2e^{-x/2}\cos\big(\frac{x\sqrt{3}}{2}-\frac{\pi}{3}\big)\big]~, \nonumber
\end{align}
also known as the third-order hyperbolic functions, and are plotted in Fig.~\ref{f1}.

\medskip

They are independent because their Wronskian $W(\phi_1,\phi_2,\phi_3)=1$, and they satisfy the following relationships
\begin{equation*}\label{nonind}
\left\{
\begin{array}{ll}
\phi_1+\phi_2+\phi_3=e^x~,\\
\phi_1^3+\phi_2^3+\phi_3^3-3\phi_1\phi_2\phi_3=1~.
\end{array}
\right.
\end{equation*}
They also fulfill the following relationships 
\begin{align}\label{ve-5}
&d\phi_1/dx=\phi_3~, \quad  d^2\phi_1/dx^2=\phi_2~,  \quad  d^3\phi_1/dx^3=\phi_1~,\nonumber\\
&d\phi_2/dx=\phi_1~, \quad d^2\phi_2/dx^2=\phi_3~, \quad  d^3\phi_2/dx^3=\phi_2~,\\
&d\phi_3/dx=\phi_2~, \quad d^2\phi_3/dx^2=\phi_1~, \quad  d^3\phi_3/dx^3=\phi_3~,\nonumber
\end{align}
which are cyclic. That is, the fourth-order derivatives have the same property as the first ones, the fifth-order derivatives as the second ones, and so forth.
In particular, they are independent solutions of the differential equation $d^3\phi/dx^3=\phi$, but also of $d^6\phi/dx^6=\phi$ and of any
$d^{3+3m}\phi/dx^{3+3m}=\phi$, with $m=0, 1,2,\dots$.

\medskip

Let
\begin{equation}\label{ve-7}
t=c\phi_1(s)+\phi_2(s)~,
\end{equation}
where $c$ is a constant. Then, $t$ fulfills
\begin{align}\label{ve-8}
&t'=c\phi_3(s)+\phi_1(s)~,\nonumber\\
&t''=c\phi_2(s)+\phi_3(s)~,\\
&t'''=t~\nonumber
\end{align}
with $\rq{}=d/ds$. Suppose $s$ is defined implicitly as a multi-valued function of $x$ as follows
\begin{equation}\label{ve-9}
x=a\left(\frac{t'}{t}\right)+b\left(\frac{t''}{t}\right)=\frac{a\phi_1(s)+bc\phi_2(s)+(ac+b)\phi_3(s)}{c\phi_1(s)+\phi_2(s)}
=\Phi(s)
\end{equation}
and denote the inverse relation by the explicit formula
\begin{equation}\label{ve-11}
s=\Phi^{-1}(x)~.
\end{equation}
Then $y=ds/dx=d\Phi^{-1}(x)/dx$ satisfies Abel's equation (\ref{Abel2}).

\medskip

Indeed, from (\ref{ve-9}), differentiating with respect to $s$ and by the usage of (\ref{ve-8}), one obtains
\begin{equation*}\label{ve-11bis}
\frac{dx}{ds}=a\bigg[\frac{t''}{t}-\left(\frac{t'}{t}\right)^2\bigg]+b\bigg[1-\left(\frac{t'}{t}\right)\left(\frac{t''}{t}\right)
\bigg]~.
\end{equation*}
The quotient $\frac{t''}{t}$ can be eliminated by means of (\ref{ve-9}), which leads to
\begin{equation}\label{ve-12}
b\frac{dx}{ds}=(ax+b^2)-(bx+a^2)\frac{t'}{t}~.
\end{equation}
Differentiating again, this time with respect to $x$ we have,
$$
-\frac{b\frac{d^2s}{dx^2}}{\left(\frac{ds}{dx}\right)^2}=a-b\left(\frac{t'}{t}\right)-(bx+a^2)\bigg[\frac{t''}{t}-
\left(\frac{t'}{t}\right)^2\bigg]\frac{ds}{dx}~.
$$
Eliminating the $t$-quotients by means of \eqref{ve-9} and \eqref{ve-12}, one can find after some reduction the following equation
$$
(bx+a^2)\frac{d^2s}{dx^2}=-2b\frac{ds}{dx}+3(ax+b^2)\left(\frac{ds}{dx}\right)^2
+(x^3-3abx-a^3-b^3)\left(\frac{ds}{dx}\right)^3~.
$$
Substituting $y=\frac{ds}{dx}$ and dividing by $(bx+a^2)$ lead to Abel's equation \eqref{Abel2}.
Thus, the solution can be written explicitly as
\begin{equation}\label{ve-15a}
\frac{b}{y}=(ax+b^2)-(bx+a^2)\bigg[\frac{c \phi_3(\Phi^{-1}(x))+\phi_1(\Phi^{-1}(x))}{c \phi_1(\Phi^{-1}(x))+\phi_2(\Phi^{-1}(x))}\bigg]~,
\end{equation}
with $c$ an arbitrary constant.

The other two solutions can be obtained by cyclically replacing the $\phi$'s in the $t$ function
\begin{equation}\label{ve-16}
t_2=c\phi_2(s)+\phi_3(s)~, \qquad t_3=c\phi_3(s)+\phi_1(s)~.
\end{equation}

\section{Vein's nonlinear oscillator}\label{sec4}

\begin{figure} [htb!]
 \centering
{\includegraphics{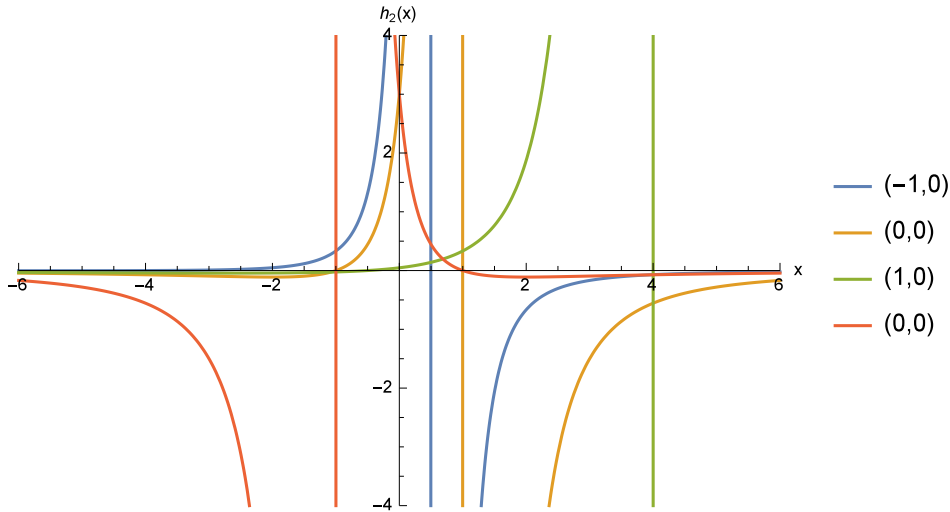}}
 \caption{The friction function $h_2(x)$ of Vein's oscillator at the fixed points $(-1,0),(0,0),(1,0)$, and $(0,0)$.}
 \label{f2}
\end{figure}

\begin{figure} [htb!]
 \centering
{\includegraphics{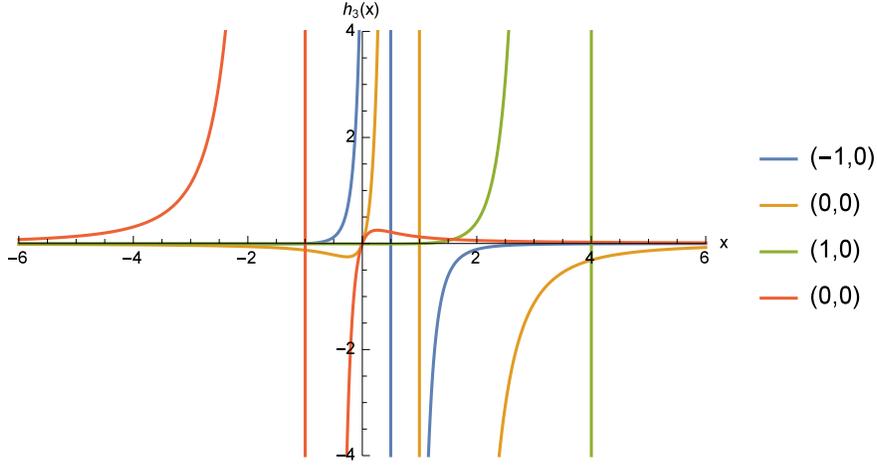}}
 \caption{The nonlinearity function $h_3(x)$ of Vein's oscillator at the fixed points $(-1,0),(0,0),(1,0)$, and $(0,0)$.}
 \label{f3}
\end{figure}

To see what kind of oscillator corresponds to Vein's Abel equation, we will proceed as follows. First, we will eliminate the linear term  from \eqref{Abel2} using
\begin{equation}
y=z e^{-2 b\int{\frac{dx}{bx+a^2}} }=\frac{z}{(bx+a^2)^2},
\end{equation}
to obtain

\begin{equation}
\frac{dz}{dx}=h_2(x)z^2+h_3(x)z^3,
\end{equation}
where
\begin{equation}\label{zee}
\begin{aligned}
	 h_2(x) &=\frac{3(ax+b^2)}{(bx+a^2)^3}~,\\
	  h_3(x) &=\frac{x^3-3abx-a^3-b^3}{(bx+a^2)^5}~.\\
	 \end{aligned}
	 \end{equation}
Then using lemma 1 with $x_{\zeta}=v(x(\zeta))=\frac {1 }{z(x(\zeta))}$, we obtain
\begin{equation}\label{x1p}
x_{\zeta \zeta}+h_2(x) x_{\zeta}  +h_3(x)=0~.
\end{equation}
which represents a nonlinear oscillator with rational friction and nonlinearity that are plotted in Figures~\ref{f2}
and~\ref{f3}, respectively, for the following values (i) $a=1, \, b=-2$; (ii) $a=1, \, b=-1$; (iii) $a=2, \, b=-1$; (iv) $a=-1, \, b=1$.

Writing equation \eqref{x1p} as a dynamical system
\begin{eqnarray}\label{l2}
\begin{array}{ll}
&x_{\zeta}=v=M(x,v)\\
&v_{\zeta}=-h_2(x)v-h_3(x)=N(x,v)~
\end{array}
\end{eqnarray}
and because \eqref{l2} cannot be put in the form of \eqref{d-9c}, because $\frac{h_3(x)}{h_2(x)}\neq const.$, we will be using instead the standard methods of phase-plane analysis and use the linear approximation at the equilibrium points of \eqref{l2} to classify them instead of solving the equation by finding the potential $\Psi$.

 Also, notice that the system is non-Hamiltonian because there is no potential $\Psi(x,v)=const.$ such that
 \begin{eqnarray}\label{s2}
\begin{array}{ll}
&\frac{\partial \Psi}{\partial  v}=v\\
&\frac{\partial \Psi}{\partial x}=h_2(x)v+h_3(x).~
\end{array}
\end{eqnarray}

   The Jacobian matrix of (\ref{l2}) is

\begin{eqnarray}\label{17}
J=\left[\begin{array}{cc}
\frac{\partial M}{\partial x}& \frac{\partial M}{\partial v}\\
\frac{\partial N}{\partial x}& \frac{\partial N}{\partial v}\\
\end{array}\right]=\left[\begin{array}{cc}
0& 1\\
-\frac{d h_2}{dx}v -\frac{d h_3}{dx}& -h_2(x)\\
\end{array}\right]~.
\end{eqnarray}
We have three equilibrium points of the system \eqref{l2} from which one is real $(x_0,y_0)=(a+b,0)$, while the other two are a pair of complex conjugates, $(x_{1,2},y_0)=(\alpha\pm i \beta,0)$, with $\alpha=-\frac{a+b}{2}$ and $\beta=-\frac{\sqrt 3 (a-b)}{2}$. The particular case $a=b$ which Vein considered, will reduce the three fixed points to the case of one real fixed point $(2a,0)$.
The characteristic polynomial of the Jacobian matrix is
\begin{equation} \label{charac1}
|J-\lambda I_2|=\lambda^2-\delta_1 \lambda +\delta_2=0~,
\end{equation}
while the discriminant is $\Delta=\delta_1^2-4\delta_2$. By evaluating the Jacobian at the real fixed point $\tilde J=J|_{(C,0)}$, where $C$ is either real or complex, we obtain

\begin{eqnarray}\label{18}
\tilde J=\left[\begin{array}{cc}
0& 1\\
-\frac{d h_3}{dx}(C)& -h_2(C)\\
\end{array}\right]~.
\end{eqnarray}
from which
\begin{equation} \label{charac2}
|\tilde J-\lambda I_2|=\lambda^2+h_2(C)\lambda +\frac{d h_3}{dx}(C)=0~.
\end{equation}
{\em 1. Real fixed point:}
In this case we have
\begin{eqnarray}\label{l3}
\begin{array}{ll}
&\delta_{1r}=-\frac{3}{(a^2+ab+b^2)^2}<0~,\\
&\delta_{2r}=\frac{3}{(a^2+ab+b^2)^4}>0~,\\
&\Delta_r=-\delta_{2r}<0~.
\end{array}
\end{eqnarray}
then, the real fixed point $(a+b,0)$  living on the $x$ axis is a stable spiral. The eigenvalues of the Jacobian matrix are  complex conjugate pairs  $\lambda_{1,2}=\frac 1 2 (\delta_1 \pm i \sqrt \delta_2)=\frac{\sqrt 3}{2(a^2+ab +b^2)^2}(-\sqrt 3 \pm i)$. By calculating the eigenvectors, the linearized solution around the fixed point is
\begin{eqnarray}\label{108}
\left[\begin{array}{cc}
\tilde x(\zeta)\\
\tilde v(\zeta)\\
\end{array}\right]=e^{-\frac{3 \zeta}{2(a^2+ab+b^2)^2}}\left[\begin{array}{cc}
\\c_1 \cos \frac{\sqrt 3 \zeta}{2(a^2+ab+b^2)^2}+c_2 \sin \frac{\sqrt 3 \zeta}{2(a^2+ab+b^2)^2}\\
\\-\frac{\sqrt 3 }{(a^2+ab+b^2)^2}\left\{c_1 \sin  \left (\frac{\sqrt 3 \zeta}{2(a^2+ab+b^2)^2}+\frac \pi 3\right)-c_2 \cos \left( \frac{\sqrt 3 \zeta}{2(a^2+ab+b^2)^2}-\frac \pi 3 \right)\right\}\\
\end{array}\right]~,
\end{eqnarray}
where $c_1$ and $c_2$ are arbitrary constants.

\medskip

{\em 2. Complex conjugated fixed points:}
In this case the coefficients of the characteristic polynomial $\delta_1,\delta_2$ are complex, and their values are shown in Table \ref{tab1} and given in (\ref{lm3}).

\begin{eqnarray}\label{lm3}
\begin{array}{ll}
&\delta_{1c}=\frac{3}{2(a^3-b^3)^2}[(a^2-2ab-2b^2)+i \sqrt 3 a (a+b)]~,\\
&\delta_{2c}=\frac{3}{2(a^3-b^3)^4}[(-a^4-8 a^3 b-6 a^2 b^2+4ab^3+2b^4)+i \sqrt 3a(a^3-6ab^2-4b^3)]~,\\
&\Delta_c=-\delta_{2c}~.
\end{array}
\end{eqnarray}

Because both coefficients $\delta_{1c}$ and $\delta_{2c}$ are not real nothing can be said about these complex fixed points.

\begin{figure} [htb!]
{\includegraphics[width=8.75cm]{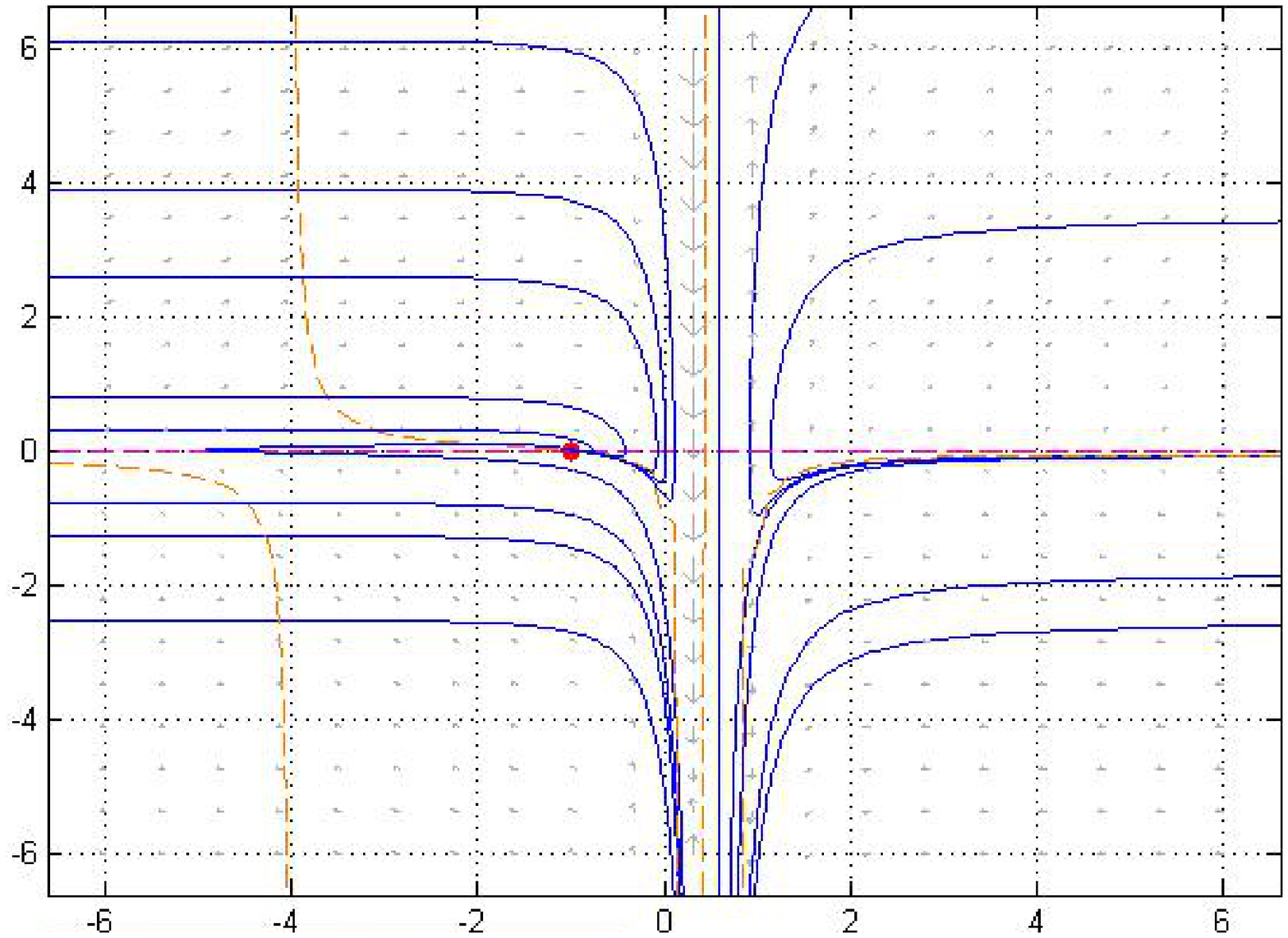}}
{\includegraphics[width=8.75cm]{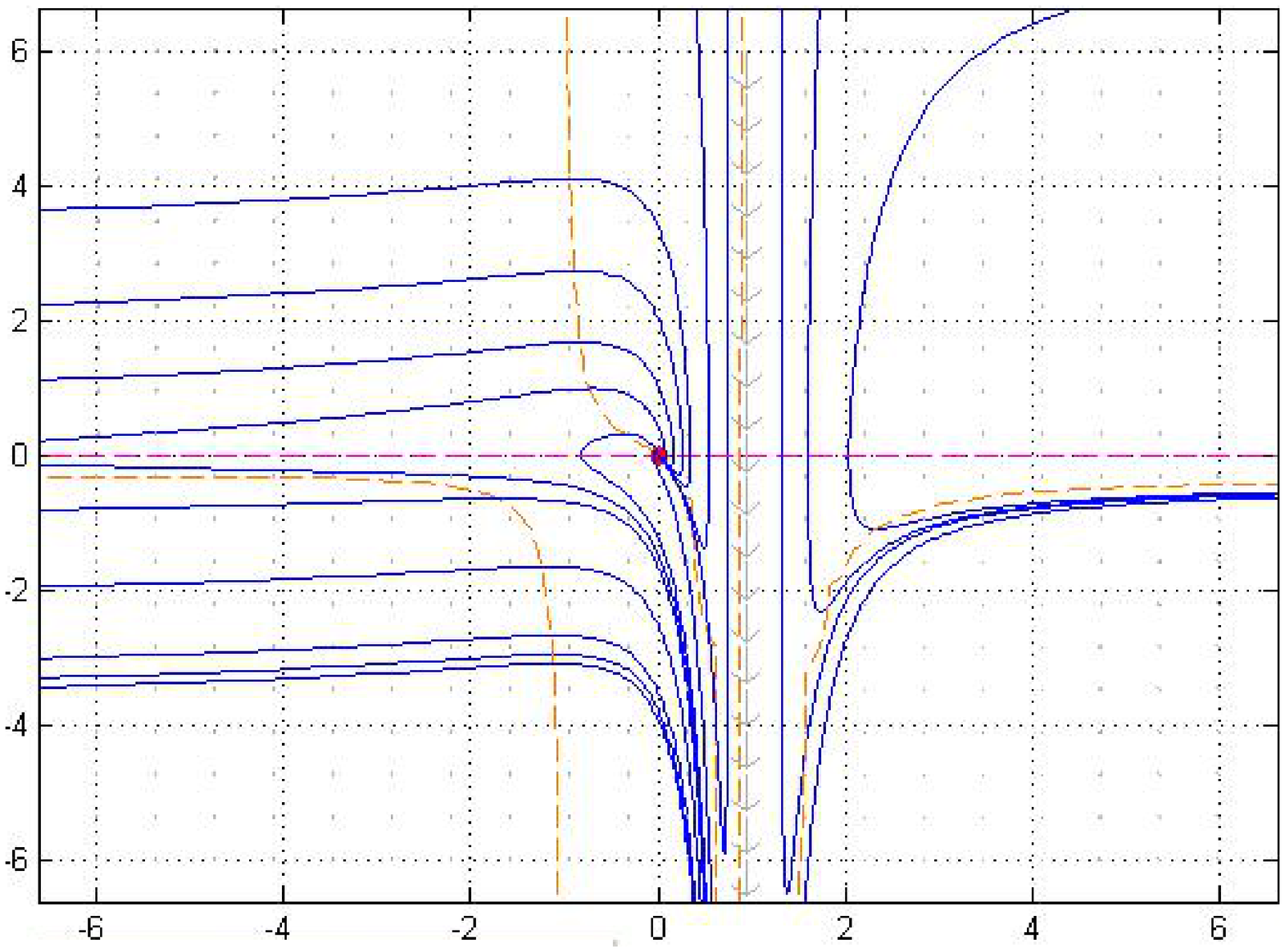}}
{\includegraphics[width=8.75cm]{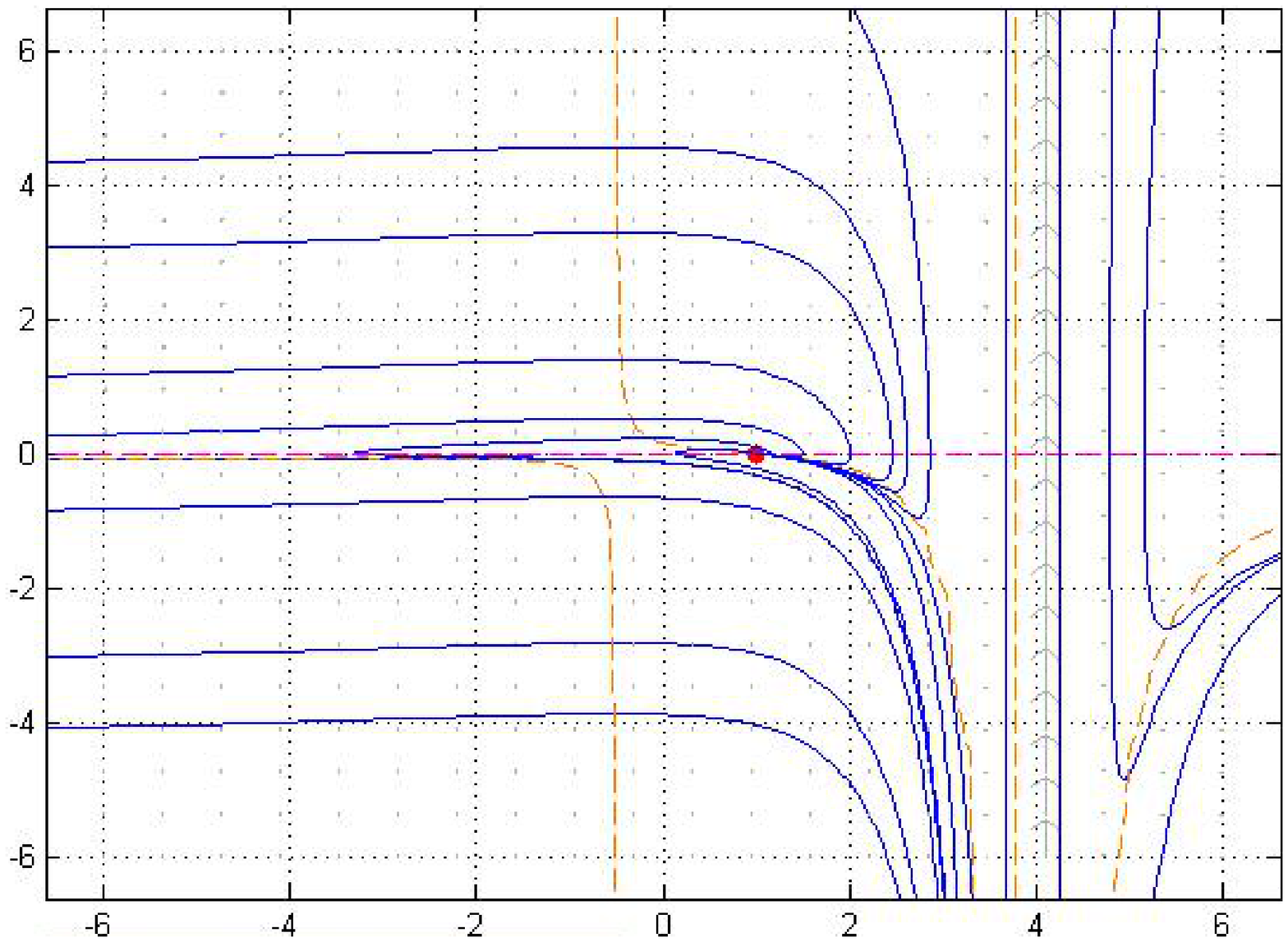}}
{\includegraphics[width=8.75cm]{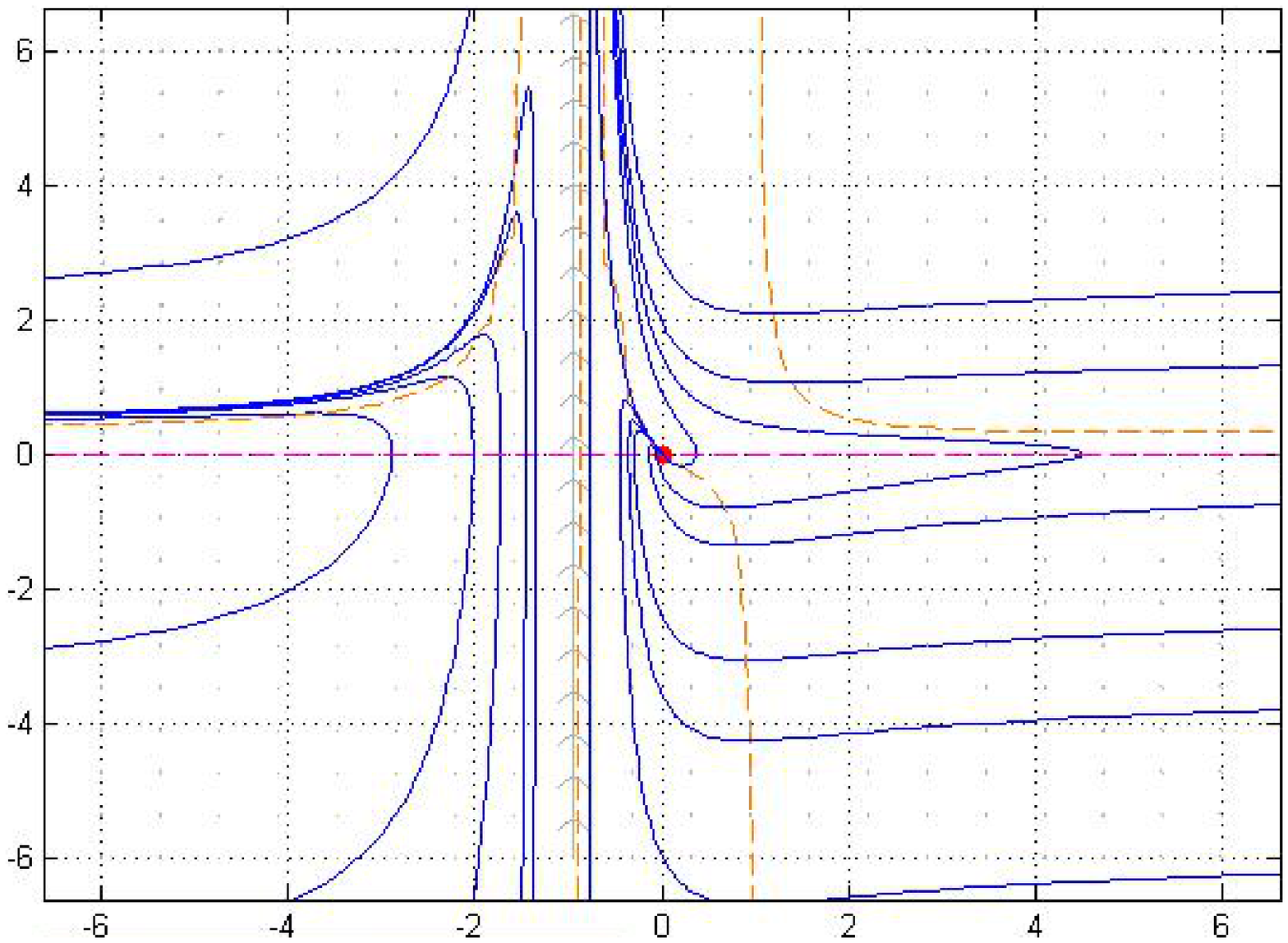}}
 \caption{Phase-plane portrait at the fixed points for $a=1, \, b=-2$ (top left); $a=1, \, b=-1$ (top right);
  $a=2, \, b=-1$ (bottom left); $a=-1, \, b=1$ (bottom right).}
 \label{f4}
\end{figure}

\bigskip
The phase-plane portrait of the system \eqref{l2} is shown in Fig. \ref{f4}. The real fixed point is identified by the red dot, and the isoclines by the
dotted curves.

\begin{table}
\begin{center}
\begin{tabular}{|c|c|c|c|c|c|}
\hline Fixed Points & $\delta_1$ & $\delta_2$ & $\Delta$ & Type \\
\hline
\hline
$(x_0,y_0)=(a+b,0)$ & $\delta_{1r}$ & $ \delta_{2r}$ & $-\delta_{2r}$ &  stable spiral \\
\hline
$(x_{1,2},y_0)=(-\frac 1 2 [(a+b)\pm i \sqrt 3 (a-b)],0)$  & $\delta_{1c}$ & $\delta_{2c}$ & $-\delta_{2c}$ & no conclusion  \\
\hline
\end{tabular}
\end{center}
\caption{Classification of equilibrium points of (\ref{l2}).}
\label{tab1}
\end{table}



%

\section{Conclusion}

This paper recalls several integrability properties of Abel equations together with some simple consequences, followed by a discussion of Vein's Abel equation related to third-order hyperbolic functions. The corresponding nonlinear oscillator system is introduced here, and its dynamical systems analysis is provided. Finally, it is worth noticing that the nonlinear oscillators corresponding to Abel equations are in fact a large class of generalized Li\'enard equations of the form $\ddot{x}+f(\dot{x})g(x)+h(x)=0$, where $g$, and $h$ are arbitrary functions and $f(\dot{x})=\alpha_1 \dot{x}^3+\alpha_2 \dot{x}^2+\alpha_3 \dot{x}+\alpha_4$, where $\alpha_i$ are constants \cite{hl} and, as such, have many applications in physics, biology, and engineering \cite{m10,nay95}.

\section{Acknowledgment}
The first author wishes to acknowledge support from Hochschule M\"{u}nchen while on leave from Embry-Riddle Aeronautical University in Daytona Beach, Florida. We also wish to thank the referees for their suggestions that led to significant improvements of this work.

\end{document}